\documentclass[12pt]{amsart}
\usepackage{amsmath,amssymb,enumerate}
\allowdisplaybreaks
\newcommand{\IB}{{\mathbb B}}

\newcommand{\IN}{{\mathbb N}}

\newcommand{\cH}{{\mathcal H}}

\newcommand{\cN}{{\mathcal N}}
\newcommand{\cU}{{\mathcal U}}

\newcommand{\ip}[1]{\mathopen{\langle}#1\mathclose{\rangle}}

\DeclareMathOperator{\Tr}{Tr}

\newtheorem*{thm}{Theorem}
\newtheorem*{lem}{Lemma}
\theoremstyle{definition}

\begin{document}
\title{A remark on amenable von Neumann subalgebras in a tracial free product}
\author{Narutaka Ozawa}
\address{RIMS, Kyoto University, \mbox{606-8502}, Japan}
\thanks{Partially supported by JSPS KAKENHI Grant Number 26400114}
\renewcommand{\subjclassname}{\textup{2010} Mathematics Subject Classification}

\begin{abstract}
Let $M=M_1*M_2$ be a nontrivial tracial free product of 
finite von Neumann algebras. We prove that any amenable 
subalgebra of $M$ that has a diffuse intersection with $M_1$ 
is in fact contained in $M_1$. 
\end{abstract}
\maketitle

Amenable von Neumann subalgebras have provided 
key tools in the study of the structure of 
the ambient von Neumann algebras.
Here we are interested in \emph{maximal} amenable subalgebras. 
The first concrete example of a maximal amenable subalgebra 
of a finite von Neumann algebra 
was found by S. Popa in \cite{popa:max}, where it is proved that the 
maximal abelian subalgebra generated by one of the generators in 
the free group factor is in fact maximal amenable. Since then most of 
the results on maximal amenability have been obtained via Popa's 
asymptotic orthogonality method. 
Recently however, R. Boutonnet and 
A. Carderi introduced in \cite{bc} a completely new way of 
proving maximal amenability via the study of central states. 
We adapt their method and prove the following. 
C.~Houdayer has shown it under more general circumstances (Theorem 4.1 in \cite{houdayer}). 

\begin{thm}
Let $M=M_1*M_2$ be a nontrivial tracial free product of 
finite von Neumann algebras $M_1$ and $M_2$, and let 
$A\subset M$ be an amenable von Neumann subalgebra 
such that $A\cap M_1$ is diffuse. 
Then, $A\subset M_1$. In particular, if $M_1$ is diffuse and 
amenable, then it is maximal amenable in $M$. 
\end{thm}

See \cite{bc,houdayer,leary} and the references therein for more information 
about maximal amenable subalgebras. 
The result should be also compared with another work 
of S. Popa (\cite{popa:orth}) that if $B\subset M_1$ is a diffuse subalgebra in the above setting, 
then the normalizer $\cN_M(B)$ of $B$ inside $M$ is contained in $M_1$. 
Here note that if $B$ is amenable in addition, then the von Neumann 
algebra $\ip{B,u}$ generated by $B$ and an element $u\in\cN_M(B)$ 
is amenable. We also mention the fascinating conjecture of J.~Peterson that 
every diffuse amenable subalgebra in the free group factor should be  
contained in a \emph{unique} maximal amenable subalgebra. 

Recall that a finite von Neumann algebra $A$ in $\IB(\cH)$ 
is \emph{amenable} if and only if there is a state $\varphi$ on $\IB(\cH)$ 
which is \emph{$A$-central}: $\varphi(ax)=\varphi(xa)$ for 
all $a\in A$ and $x\in \IB(\cH)$.

\begin{lem}
Let $A\subset \IB(\cH)$ be a finite amenable von Neumann subalgebra 
and $x\in\IB(\cH)$ be such that the norm closed convex hull of 
$\{ (u\otimes\bar{u}) (x\otimes\bar{x})  (u\otimes\bar{u})^* : u\in\cU(A)\}$ 
in $\IB(\cH\otimes\bar{\cH})$ contains $0$.
Then, $\varphi(x^*A)=\{0\}$ for any $A$-central state $\varphi$. 
\end{lem}
The condition on $x$ is reminiscent of singularity in \cite{bc}. 
Indeed, it can be shown that a subgroup $\Lambda\le\Gamma$ is \emph{singular} 
(Definition 1.2 in \cite{bc}) 
if and only if for every $g\in\Gamma\setminus\Lambda$ the closed convex hull of 
$\{ \lambda(tgt^{-1}) : t\in\Lambda\}$ in $\IB(\ell_2\Gamma)$ contains $0$. 

\begin{proof}[Proof of Lemma]
Let $\varphi$ be an $A$-central state on $\IB(\cH)$ and approximate it 
by a net $(S_i)_i$ of positive norm one trace class operators. 
We may assume that $\|[S_i,a]\|_1\to0$ for all $a\in A$, by the standard trick of Day.
Thus, for every unitary element $a\in A$, one has 
\[
|\varphi(x^*a)| = \lim_i|\Tr(S_ix^*a)| 
 = \lim_i|\Tr(S_i^{1/4}x^*S_i^{1/4}aS_i^{1/2})| 
 \le \limsup_i|\Tr(S_i^{1/2}x^*S_i^{1/2}x)|^{1/2}.
\]
Here we have used the Kittaneh--Powers--St{\o}rmer inequality (Corollary~2 in \cite{kittaneh}) 
\[
\| aS_i^{1/4}a^* - S_i^{1/4}\|_4^4
\le\| aS_i^{1/2}a^* - S_i^{1/2}\|_2^2
\le\| aS_ia^*-S_i\|_1\to0.
\]
For any finite (multi-)subset $F$ of unitary elements in $A$, one has 
\begin{align*}
\limsup_i |\Tr(S_i^{1/2}x^*S_i^{1/2}x)|
 &= \limsup_i | \frac{1}{|F|}\sum_{u\in F}\Tr(u^*S_i^{1/2}ux^*u^*S_i^{1/2}ux) | \\
 &\le \limsup_i \| \frac{1}{|F|}\sum_{u\in F} uxu^*S_i^{1/2}ux^*u^*\|_2 \\
 &\le \|\frac{1}{|F|}\sum_{u\in F}(uxu^*)\otimes(\bar{u}\bar{x}\bar{u}^*) \|_{\IB(\cH\otimes\bar{\cH})}.
\end{align*}
Here we have used the identification of 
the Hilbert--Schmidt class operators and the vectors in $\cH\otimes\bar{\cH}$.
Now, the assumption implies that $\varphi(x^*a)=0$.
\end{proof}

\begin{proof}[Proof of Theorem]
Since $A$ is an amenable subalgebra of $M$, the tracial state on $M$ 
extends to an $A$-central state $\varphi$ on $\IB(L^2M)$.
It suffices to show $\varphi(x^*A)=\{0\}$ for every $x$ in $M\ominus M_1$. 
We may assume that $x$ is of the form $v_1\cdots v_l$ 
for some trace zero unitary elements $v_j\in M_{i(j)}$ with $i(j)\neq i(j+1)$. 
Since $A\cap M_1$ is diffuse, it contains a unitary element $u$ 
whose nonzero powers are all trace zero. 
Thus $a:=x^*ux$ and $b:=u$ are free Haar unitary elements, 
and they together generate a copy of the free group factor $LF_2$ of rank two. 
This implies that $\{ a^kb^{-k} : k\in\IN\}$ is also a free family.
It follows that 
\[
\|\frac{1}{m}\sum_{k=1}^m (x^*u^kxu^{-k}) \otimes(\bar{x}^*\bar{u}^k\bar{x}\bar{u}^{-k})\| 
=\|\frac{1}{m}\sum_{k=1}^m a^kb^{-k}\|_{LF_2}
=\frac{2\sqrt{m-1}}{m}\to0
\]
by the Akemann--Ostrand formula (\cite{ao}). Hence $\varphi(x^*A)=\{0\}$ by the above lemma.
\end{proof}

\end{document}